\documentclass[reqno]{amsart}
\usepackage{amssymb, latexsym}

\usepackage{xspace}
\usepackage[bookmarksnumbered,colorlinks]{hyperref}
\usepackage{graphics}
\def\bb#1\eb{\textcolor{blue}
{#1}} %
\def\br#1\er{\textcolor{red}
{#1}} %
\def\bv#1\ev{\textcolor{green}
{#1}} %
\def\bc#1\ec{\textcolor{cyan}
{#1}} %

\usepackage{graphics}

\def\bal#1\eal{\begin{align}#1\end{align}}                      %
\def\baln#1\ealn{\begin{align*}#1\end{align*}}          
\def\bml#1\eml{\begin{multline}#1\end{multline}}        %
\def\bmln#1\emln{\begin{multline*}#1\end{multline*}}  %
\def\bga#1\ega{\begin{gather}#1\end{gather}}
\def\bgan#1\egan{\begin{gather*}#1\end{gather*}}

\makeatletter

\@addtoreset{equation}{section} \makeatother
\newtheorem{theorem}{Theorem}[section]

\newtheorem{corollary}[theorem]{Corollary}

\newtheorem{definition}[theorem]{Definition}

\newtheorem{remark}[theorem]{Remark}

\newcommand{\dimo}{\noindent{\em Proof.}$\ $}

\newcommand{\N}{\ensuremath{\mathbb{N}}\xspace}     
\newcommand{\R}{\ensuremath{\mathbb{R}}\xspace}     
\newcommand{\Z}{\ensuremath{\mathbb{Z}}\xspace}     
\newcommand{\eps}{\varepsilon}                      

\title[Perturbed asymptotically linear problems]{Perturbed asymptotically linear problems}
\thanks{Partially supported by M.I.U.R. Research project 
PRIN2009 ``Metodi Variazionali  e Topologici nello Studio di Fenomeni Nonlineari''}

\author[R. Bartolo]{R. Bartolo}
\address{Rossella Bartolo\hfill\break\indent
Dipartimento di Ingegneria Meccanica e Gestionale\hfill\break\indent
Politecnico di Bari\hfill\break\indent
Via E. Orabona 4, 70125 Bari\hfill\break\indent
Italy}
\email{r.bartolo@poliba.it}

\author[A.M. Candela]{A.M. Candela}

\author[A. Salvatore]{A. Salvatore}
\address{Anna Maria Candela and Addolorata Salvatore\hfill\break\indent
Dipartimento di Matematica\hfill\break\indent
Universit\`a degli Studi di Bari ``Aldo Moro''\hfill\break\indent 
Via  Orabona 4, 70125 Bari\hfill\break\indent
Italy}
\email{candela@dm.uniba.it}
\email{salvator@dm.uniba.it}



\date{}

\begin{document}

\begin{abstract}
The aim of this paper is investigating the existence of solutions of the semilinear elliptic problem
\begin{equation}
\left\{
\begin{array}{ll}
\displaystyle{-\Delta u\ =\ p(x, u) + \varepsilon g(x, u)} 
& \mbox{ in }  \Omega,\\
 \displaystyle{u=0}  & \mbox{ on }  \partial\Omega,\\
 \end{array}
\right.
\end{equation}
where $\Omega$ is an open bounded domain of $\R^N$, $\varepsilon\in\R$, $p$ is subcritical 
and asymptotically linear at infinity and $g$ is just a continuous function. Even when this 
problem has not a variational structure on $H^1_0(\Omega)$, suitable procedures and estimates 
allow us to prove that the number of distinct crtitical levels of the functional associated to 
the unperturbed problem is ``stable'' under small perturbations, in particular obtaining 
multiplicity results if $p$ is odd, both in the non--resonant and in the resonant case.
\end{abstract}

\maketitle
\section{Introduction}

Let us consider the following semilinear elliptic problem
\begin{equation}\label{problema}
\left\{
\begin{array}{ll}
\displaystyle{-\Delta u\ =\ p(x, u) + \varepsilon g(x, u)} 
& \mbox{ in }  \Omega,\\
 \displaystyle{u=0}  & \mbox{ on }  \partial\Omega,\\
 \end{array}
\right.
\end{equation}
where $\Omega$ is an open bounded domain of $\R^N$ ($N\geq 3$)
with smooth boundary $\partial\Omega$, $\eps \in\R$ 
and $p,g$ are given real functions on $\Omega\times \R$. 

If $\eps =0$ problem \eqref{problema} has 
been widely investigated when $p(x,\cdot)$ is asymptotically linear and possibly odd
 (see \cite{az,bbf} and references therein). 

On the other hand, if $\eps\not=0$, let us
consider a perturbation term 
$g$ just continuous, without assumptions on its growth. In \cite{ll} Li and Liu state  the existence 
of multiple solutions of \eqref{problema} 
when $p(x,\cdot)$ is odd, superlinear at infinity, but subcritical
(see also \cite{dr} and \cite{ll1} for related results). 
Furthermore, in \cite[Theorem 1.6]{ll1} a multiplicity result is stated if $p(x,\cdot)$ 
is asymptotically linear at infinity and odd when $g(x,\cdot)$ is odd, too (see also 
Theorem \ref{dis} below). Here, our aim is going further in this direction and proving 
that, even when \eqref{problema} has not a variational structure on $H_0^1(\Omega)$,
the number of distinct critical levels of the functional associated to the unperturbed 
problem is ``stable'' under small perturbations also in lack of symmetry.

Throughout this paper suppose that 
there exist $\lambda \in \R$ and $f:\Omega\times\R\rightarrow \R$ such that
\begin{equation}\label{asin}
p(x,t)\ =\ \lambda t + f(x,t)\quad \hbox{for all $t\in\R$ and a.e. $x\in \Omega$.}
\end{equation}
Denote respectively by $\sigma(-\Delta)$ and by
$0<\lambda_1<\lambda_2<\ldots<\lambda_k<\ldots$ the spectrum and the distinct 
eigenvalues of $-\Delta$ in $H^1_0(\Omega)$. 

Let us introduce the following conditions:
\begin{itemize}
\item[$(H_0)$] 
$f$ is a Carath\'eodory function (i.e. $f(\cdot,t)$ is measurable in $\Omega$
 for all $t\in\R$ and $f(x,\cdot)$ is continuous in $\R$ for a.e. $x\in \Omega$) 
 and
\begin{equation}\label{sup}
\sup_{|t|\leq r}|f(\cdot,t)|\in L^\infty(\Omega) \quad \hbox{ for all }\quad r>0;
\end{equation}
\item[$(H_1)$] there exists
\[
\lim_{|t|\to + \infty}\frac{f(x,t)}{t}\ =\ 0\quad \hbox{uniformly with respect to a.e. $x \in \Omega$;}
\]
\item[$(H_2)$]
there exists
\[
\lim_{t\rightarrow 0}\frac{f(x,t)}{t} \ =:\ f'(0) \in\R \quad \hbox{uniformly with respect to a.e. 
$x \in \Omega$;}\]
\item[$(H_3)$] $\lambda\not\in\sigma(-\Delta)$;
\item[$(H_4)$] there exist two integers $h$, $k \ge 1$ such that
\[
\min\{f'(0) + \lambda,\lambda\}\ <\ \lambda_h\ < \ \lambda_k\ <\max\{f'(0) + \lambda,\lambda\}; 
\]
\item[$(H_5)$] $f(x,\cdot)$ is odd for a.e. $x\in\Omega$.
\end{itemize}

\vspace{0,5cm}
By \eqref{asin} problem \eqref{problema} becomes
\[
\leqno{(P_{\varepsilon})}\qquad\qquad
\left\{
\begin{array}{lll}
\displaystyle{ -\Delta u -\lambda u=  f(x, u) + \varepsilon g(x, u)  } & & \mbox{ in }  \Omega,\\
 \displaystyle{u=0} & & \mbox{ on }  \partial\Omega.\\
 \end{array}
\right.
\]

If $(H_0)$--$(H_1)$ hold and $\varepsilon =0$, it is well known that the 
solutions of problem $(P_0)$ are the critical points of the $C^1$ functional
\begin{equation}
\label{1}
I(u)=\frac 12\int_\Omega |\nabla u|^2\; {\rm d}x - \frac\lambda 2 
\int_\Omega u^2 \; {\rm d}x - \int_\Omega F(x,u) \;{\rm d}x \quad \hbox{ on } \quad H^1_0(\Omega),
\end{equation}
with $F(x,t)=\int_0^t f(x,s) \; {\rm d}s$.
By using the pseudo--index theory introduced by Benci in \cite{b}, it can 
be proven that $I$ has at least as many distinct critical points as the 
number of the eigenvalues of $-\Delta$ between $f'(0)+ \lambda$ and $\lambda$, 
even if some of the corresponding critical levels may be the same
(cf. \cite{bbf} and, here, Section \ref{s21}). 

In our first theorem we give an existence result for $(P_{\varepsilon})$ and, 
in the easier case in which both $f$ and the pertubation function $g$ are odd, 
we obtain also multiplicity of solutions: the crucial point is that, in this case, 
near critical values of $I$ there are critical values of odd perturbations of $I$ (see Theorem \ref{near}).

\begin{theorem}\label{maintheorem}
Assume that $g\in C(\overline\Omega\times\R, \R)$ and $(H_0)$, $(H_1)$ and $(H_3)$ hold. 
Then, there exists $\bar\varepsilon >0$ such that problem $(P_{\varepsilon})$ has 
at least one solution for all $|\varepsilon|\leq\bar\varepsilon$. \\
Moreover, if $g$ is odd, under the further assumptions $(H_2), (H_4)$ and $(H_5)$, 
denoting by $\bar m$, with $1\leq \bar m\leq \dim (M_h\oplus\ldots \oplus M_k)$, 
the number of the distinct mini--max critical levels of the unperturbed functional 
$I$ in \eqref{1}, then  problem $(P_{\varepsilon})$ has at least $\bar m$ distinct 
pairs of solutions for all $|\varepsilon|\leq\bar\varepsilon$.
\end{theorem}

We point out that assumption $(H_3)$ can be avoided. Indeed, many authors studied 
problem $(P_0)$ in the so called ``resonant case'', but under suitable stronger 
hypotheses on the nonlinear term $f$ (see \cite{bbf} and references therein).
Here, we consider the following further conditions:
\begin{itemize}
\item[$(H_1')$] there exists $M>0$ such that
\[|f(x,t)|\leq M \quad \hbox{for a.e. $x \in \Omega$ and for all $t\in\R$;}\]
\item[$(H_3')$] there exists  an integer $k \ge 1$ such that $\lambda=\lambda_k\in\sigma(-\Delta)$;
\item[$(H_4')$] there exists an integer $h\ge 1$ such that
\[
\min\{f'(0) + \lambda_k,\lambda_k\}\ <\ \lambda_h\ <\max\{f'(0) + \lambda_k,\lambda_k\}; 
\]
\item[$(H_6)$] there exists
\[\lim_{|t|\to + \infty}\ F(x,t)=\ 	l\in \{\pm\infty\}\quad \hbox{uniformly with respect to a.e. $x \in \Omega$.}\]
\end{itemize}
\begin{theorem}\label{strong}
Assume that $g\in C(\overline\Omega\times\R, \R)$ and $(H_0), (H_1'), (H_3')$ and $(H_6)$ hold.  
Then there exists $\bar\varepsilon >0$ such that problem $(P_{\varepsilon})$ has at least one 
solution for all $|\varepsilon|\leq\bar\varepsilon$. \\
Moreover, if $g$ is odd, under the further assumptions $(H_2), (H_4')$ and $(H_5)$, 
the same multiplicity result of Theorem \ref{maintheorem} holds.
\end{theorem}

We point out that in \cite[Theorem 1.6]{ll1} the authors obtain a result 
sharper than our in Theorem \ref{strong}, but it works only in the non 
resonant case. Furthermore, our approach allows us to deal also with 
the case in which the symmetry of problem $(P_{\varepsilon})$ is broken 
by a non--odd perturbation term $g$. In this case, we still obtain some 
multiplicity results, but our techniques impose to work with nonlinearities 
such that the functional $I$ in \eqref{1} has only critical levels which arise 
for topological reasons, called {\em topologically relevant} (see Definition \ref{toprel}). 
Anyway, similar assumptions appear also in previous references concerning asymptotically 
linear (unperturbed) Dirichlet problems, as different approaches require sometimes that 
the critical points of the associated functional are ``nondegenerate'' 
(see e.g. \cite{cos} and references therein). 
Roughly speaking, we need to restrict a priori 
the set of the critical values of our functional $I$ in \eqref{1}, 
so that they are indeed {\em essential} ones (cf. Definition \ref{essential1}) 
and a stability result, stating that critical values of $I$ are preserved for 
small perturbations, holds (see Corollary \ref{near2}, while we refer to 
\cite[Theorem 8.10]{mw} for a related result based on Morse Theory).
In \cite{rei} Reeken obtains something stronger. Indeed, he gives a topological 
description of some critical levels $c$ of a $C^1$ functional and proves that 
the category of the set of the critical points of small perturbations of the 
functional is lower bounded from an integer representing the topological 
description of $c$ (see \cite[Theorem 5.2]{rei}).

\begin{theorem}\label{dis}
Assume that $g\in C(\overline\Omega\times\R, \R)$ and $(H_0)$--$(H_5)$ hold. 
If all the critical levels of the functional  $I$ in (\ref{1}) are topologically 
relevant, denoting by $\bar m$, with $1\leq \bar m\leq \dim (M_h\oplus\ldots \oplus M_k)$, 
the number of the distinct mini--max critical levels of $I$, then  problem 
$(P_{\varepsilon})$ has at least $\bar m$ distinct solutions for all $|\varepsilon|\leq\bar\varepsilon$.
\end{theorem}

We can extend the previous theorem to the resonant case as follows.

\begin{theorem}\label{disst}
Under the assumptions of Theorem \ref{dis}, just replacing hypotheses 
$(H_1)$, $(H_3)$, $(H_4)$ respectively by $(H_1')$, $(H_3')$, $(H_4')$ 
and adding condition $(H_6)$, the same multiplicity result of Theorem \ref{dis} holds also in the resonant case.
\end{theorem}

\begin{remark}\label{dora1}
{\em The previous results hold also if in assumption $(H_2)$ it is $f'(0)=0$. 
Indeed, in this case $(H_4)$ is not meaningful, but imposing a condition on 
the sign of $F$, a multiplicity result can still be stated (see Remarks \ref{inzero} 
and \ref{dora3}). On the other hand, our theorems can be also proven when the limit 
in $(H_2)$ is infinite, i.e. $f'(0)\in \{\pm \infty\}$ (see Remarks \ref{dora2} and \ref{dora3} for more details).}
\end{remark}

This paper is organized as follows.
In Section \ref{s22} we recall some abstract tools, 
in particular concerning the notions of pseudo--index and essential value,
in Section \ref{s21} we deal with the unperturbed problem $(P_0)$ and, 
lastly, in Section \ref{s23} we prove our main results.
\vspace{0,5cm}

\noindent
{\sl Notations.} Throughout this paper we denote by $(X,\|\cdot\|_X)$ a Banach space, 
by $(X',\|\cdot\|_{X'})$ its dual, by $J$ a $C^1$ functional on $X$ and by
\begin{itemize}
\item $J^b = \{u\in X : J(u)\leq b\}$ the sublevel of $J$ corresponding to $b\in \bar\R:=\R\cup\{\pm\infty\}$;
\item $K_c  = \{u \in X:\ J(u) = c,\ {\rm d}J(u) = 0\}$ the set of the critical points of $J$ 
in $X$ at the critical level $c \in \R$.
\end{itemize}
Furthermore, let us denote by
\begin{itemize}
\item $|\cdot|_s$ the usual norm in the Lebesgue space $L^s(\Omega)$, $1\leq s\leq +\infty$;
\item $\|\cdot\|$ the norm in $H^1_0(\Omega)$, i.e. $\|u\|=|\nabla u|_2$ for all $u\in H^1_0(\Omega)$;
\item $2^\ast= \frac{2N}{N-2}$ the critical exponent for Sobolev embeddings of $H^1_0(\Omega)$ $(N\geq 3)$;
\item  $M_{j}$ the eigenspace corresponding to the eigenvalue $\lambda_j$ of $-\Delta$ in $H^1_0(\Omega)$  
and  $u_j$ the component of $u$ in $M_j$,  for any integer $j\ge 1$ and for each $u\in H^1_0(\Omega)$;
\item
$H^-(j)\ =\ {{\displaystyle\bigoplus_{i\leq j}M_{i}}}$ and $H^+(j)=
\overline{\displaystyle\bigoplus_{i\geq j}M_{_i}}$, for any integer $j\ge 1$;
\item $C_j$ a positive real number, for any integer $j\ge 1$.
\end{itemize}

\section{Some abstract tools}\label{s22}

We recall the well--known Palais--Smale condition.
 
\begin{definition}\label{ps}
{\rm The functional $J$ satisfies the {\slshape Palais--Smale Condition at
level $c$} ($c \in \R$), briefly $(PS)_c$,
if any sequence $(u_n)_n \subset X$ such that}
\[
\lim_{n \to +\infty}J(u_n) = c\quad\mbox{ \rm and}\quad
\lim_{n \to +\infty}\|{\rm d}J(u_n)\|_{X'} = 0,
\]
{\rm converges in $X$, up to subsequences. In general, 
if $-\infty\leq a<b\leq +\infty$, $J$ satisfies $(PS)$ in $]a,b[$ if so is at each level $c\in ]a,b[$.}
\end{definition} 

Beside the existence critical point theorems, under this condition 
sharper results can be proven when one deals with symmetric functionals 
(see e.g. \cite{ar}). In \cite{bbf} multiplicity results for critical 
points of even functionals are stated and their proofs are based on the 
use of a pseudo--index theory. In order to introduce such definition, 
let us recall some notions of the index theory
for an even functional with symmetry group $\Z_2 = \{{\rm id}, -{\rm id}\}$.

Define
\[
\begin{split}
\Sigma = \Sigma(X) \ =\ \{A \subset X:\ &A \ \hbox{closed and symmetric w.r.t. the origin,}\\
&\hbox{i.e. $-u \in A$ if $u \in A$}\}\end{split}
\]
and
\[
{\mathcal H} =\{h\in C(X,X): h \mbox{ odd} \}.
\]
Taking $A \in \Sigma$, $A\ne\emptyset$, the genus of $A$ is
\[
\gamma(A) \ =\ \inf\{k \in \N^*:\  \exists \psi \in C(A,\R^k\setminus\{0\})\ \hbox{s.t.}\
\psi(-u) = - \psi(u) \text{ for all } u\in A\},
\]
if such an infimum exists, otherwise $\gamma(A) = +\infty$. Assume $\gamma(\emptyset) = 0$.

The index theory $(\Sigma,{\mathcal H}, \gamma)$ related to $\Z_2$ is also 
called {\em genus} (see \cite{Kr} and for more details
\cite[Section II.5]{s}).

According to \cite{b}, the pseudo--index related to the genus, an 
even functional $J:X\rightarrow \R$ and $S\in\Sigma$ is the triplet
$(S, {\mathcal H}^\ast,\gamma^\ast)$ such that
\[
{\mathcal H}^\ast =\{h\in{\mathcal H}: h \mbox{ bounded homeomorphism s.t. } h(u)=u \mbox{ if } u\not\in J^{-1}(]0,+\infty[) \}
\]
and
\[
\gamma^\ast(A) =\min _{h\in {\mathcal H}^\ast}\gamma(h(A)\cap S)  \;\; \text{ for all } A\in\Sigma.
\]

A {\em mini--max} theorem was stated in \cite[Theorem 2.9]{bbf} 
under the weaker Cerami's variant of the Palais--Smale condition 
(cf. e.g. \cite[Definition 1.1]{bbf}), while here we recall it when just the $(PS)$ condition holds.

\begin{theorem}\label{group}
Let $H$ be a real Hilbert space, $J\in C^1(H,\R)$ an even functional, $(\Sigma,\mathcal H,\gamma)$ 
the genus theory on $H$, $S\in\Sigma$, $(S, {\mathcal H}^\ast,\gamma^\ast)$ the pseudo--index 
theory related to the genus, $J$ and $S$. Consider $a,b, c_0, c_\infty\in\bar\R$, 
$-\infty\leq a<c_0<c_\infty<b\leq +\infty$.
Assume that:
\begin{list}{(\roman{enumi})}{\usecounter{enumi}\labelwidth  5em
\itemsep 0pt \parsep 0pt}
\item the
functional $J$ satisfies $(PS)$ in $]a,b[$;
\item $S\subset J^{-1}([c_0,+\infty[)$; 
\item there exist an integer $\bar k\geq 1$ and $\bar A\in\Sigma$ such that 
$\bar A\subset J^{c_\infty}$ and $\gamma^\ast(\bar A)\geq \bar k$.
\end{list}
Then the numbers
\begin{equation}\label{value}
c_i=\inf_{A\in \Sigma_i}\sup_{u\in A}J(u), \quad \quad i\in \{1,\ldots, \bar k\},
\end{equation}
with $\Sigma_i= \{A\in \Sigma: \gamma^\ast(A)\geq i\}$,
are critical values for $J$ and
\[
c_0\leq c_1\leq \ldots\leq c_{\bar k}\leq c_\infty.
\] 
Furthermore, if $c=c_i=\ldots =c_{i+r}$, with $i\geq 1$ and $i+r\leq \bar k$, then $\gamma(K_c)\geq r+1$. 
\end{theorem}

We have already pointed out that in our main theorems we may deal with problems  
$(P_\varepsilon)$ without a variational structure on $H^1_0(\Omega)$. Hence, 
following \cite{ll}, we use the auxiliary notion of {\em essen\-tial va\-lue}, 
as it is introduced in \cite{dl1} (see also \cite{dl}) in the study of 
perturbations of nonsmooth functionals. Moreover, for even functionals 
we introduce  the definition of {\em odd--essential value}, which allows 
us to obtain multiplicity results for odd perturbations.

\begin{definition}\label{essential}
{\rm Let $J: X\rightarrow \R$ be continuous (resp. $J$ even continuous) and $a,b\in\bar\R$, with $a\leq b$.
The pair 
$(J^b,J^a)$ is {\slshape trivial} (resp. {\slshape odd--trivial}) if,  
for each neighbourhood  $[\alpha',\alpha'']$ 
of $a$ and $[\beta',\beta'']$ of $b$ in $\bar \R$, there exists a 
continuous (resp. an odd continuous) map $\varphi: J^{\beta'}\times [0,1]\rightarrow J^{\beta''}$ such that}
\begin{itemize}
\item[$(i)$]
$\varphi (x,0)=x$  {\rm for each} $x\in J^{\beta'};$
\item[$(ii)$]
$\varphi (J^{\beta'}\times\{1\})\subseteq J^{\alpha''};$
\item[$(iii)$]
$\varphi(J^{\alpha'}\times [0,1])\subseteq J^{\alpha''}.$
\end{itemize}
\end{definition} 

\begin{definition}\label{essential1}
{\rm Let $J: X\rightarrow \R$ be continuous (resp. $J$ even continuous).
A real number $c$ is an {\slshape essential value}} 
{\rm (resp. an {\slshape odd--essential value}) of $J$ if for each $\varepsilon>0$ 
there exist }$a,b\in ]c-\varepsilon, c+\varepsilon[$, $a<b$, {\rm such that the pair 
$(J^b,J^a)$ is not trivial (resp. not odd--trivial).}
\end{definition} 

The following theorem states that small perturbations of a continuous functional preserve 
the essential values (see \cite[Theorem 3.1]{dl1} or also \cite[Theorem 2.6]{dl}); 
in particular this holds for the odd ones, just doing some small changes in the proof of \cite[Theorem 3.1]{dl1}.

\begin{theorem}\label{near}
Let $c\in\R$ be an essential value (resp. odd--essential value) of 
$J: X\rightarrow \R$ continuous (resp. $J$ even continuous). Then,
for every $\eta >0$ there exists $\delta>0$ such that every functional (resp. even functional) $G\in C(X,\R)$ with
\[
\sup\{|J(u)- G(u)|: u\in X\} <\delta
\]
admits an essential value (resp. odd--essential value) in $]c-\eta, c+ \eta[$.
\end{theorem}

Now, we focus on the setting of smooth functionals and 
recall some results which link critical and essential values,  
stating in particular that the critical values arising from 
mini--max procedures are essential, if all the involved deformations 
are of the ``same kind'' (see \cite[Theorems 3.7 and 3.9]{dl1}). 

\begin{theorem}\label{near1}
Let $c\in\R$ be an essential value of $J\in C^1(X,\R)$. If $(PS)_c$ holds, then
$c$ is a critical value of $J$.
\end{theorem}

\begin{remark}\label{example}
{\em In general the reverse implication does not hold: even if $(PS)_c$ is satisfied, 
a critical value is not necessarily an essential one (see e.g. \cite[Example 3.12]{dl1}).}
\end{remark}

\begin{theorem}\label{closed}
Let $\Gamma$ be a non empty family of non empty subsets of $X$, 
$J\in C^1(X,\R)$ and $d\in\R\cup \{-\infty\}$. Let us assume that, 
for every $C\in\Gamma$ and for every deformation $\varphi:  X\times [0,1]\longrightarrow X$ 
with $\varphi(u,t)=u$ on $J^d\times [0,1]$, it is $\overline{\varphi(C\times \{1\})}\in\Gamma$.\\
Setting 
\begin{equation}\label{val}
c=\inf_{C\in\Gamma}\sup_{u\in C} J(u),
\end{equation}
if $d<c<+\infty$, then $c$ is an essential value of $J$.
\end{theorem}

Let us point out that, if $J$ is even, the previous theorem does not apply to the 
critical values $c_i$ given by \eqref{value}; in fact, if $\varphi$ is a 
deformation as in Theorem \ref{closed}, the set $\overline{\varphi(C\times \{1\})}$, 
$C\in\Gamma$, does not necessarily belong to $\Gamma$ because $\varphi$ could be not 
odd. Thus, we cannot assert that the $c_i$'s are indeed essential values of $J$. 
Anyway, slight modifications in the proof of \cite[Theorem 3.9]{dl1} allow us to 
state the following result concerning odd--essential values.

\begin{corollary}\label{near5}
Let $\Gamma$ be a non empty family of non empty symmetric subsets of $X$, 
$J\in C^1(X,\R)$ even and $d\in\R\cup \{-\infty\}$. Let us assume that, 
for every $C\in\Gamma$ and for every odd deformation $\varphi:  X\times [0,1]\longrightarrow X$ 
with $\varphi(u,t)=u$ on $J^d\times [0,1]$, it is $\overline{\varphi(C\times \{1\})}\in\Gamma$.\\
Then, taking  $c$  as in \eqref{val}, if $d<c<+\infty$, we have that $c$ is an odd--essential value of $J$.
\end{corollary}

As we have observed in the Introduction, we also deal with perturbation from symmetry problems, 
thus we restrict ourselves to consider a subset of the critical values of $I$ in \eqref{1}, so 
that a stability result does still hold (we also refer to \cite{hz} for a related result 
concerning changing--sign solutions of some elliptic equations).

Actually, in Corollary \ref{near2} below we consider just the 
preservation for small perturbations of some critical levels of 
a smooth functional satisfying the $(PS)$ condition.
>From the Deformation Lemma, if $J\in C^1(X,\R)$ satisfies $(PS)_c$ 
and $c$ is not a critical value, then for any $\bar\eta>0$ there 
exists $\eta\in ]0,\bar\eta[$ and $\varphi\in C(X\times [0,1],X)$, 
which is odd if $J$ is even, such that $\varphi(u,1)=u$ if 
$J(u)\not\in [c-\bar\eta, c+ \bar\eta]$ and $\varphi(J^{c+\eta},1)\subset J^{c-\eta}$. 
In some sense now we require that also the other implication is true, that is if 
$J^{c-\eta}$ is a strong deformation retract of $J^{c+\eta}$, then $c$ is not critical.

Indeed, starting from Definitions \ref{essential} and \ref{essential1}, we give the following definition.

\begin{definition}\label{toprel}{\rm 
A critical level $c$ of a  functional $J\in C^1(X,\R)$ is {\em topologically relevant} if it is an essential one.}
\end{definition}

Summing up, we work with a special class of critical levels: those which are also essential.

Thus, according to Definition \ref{toprel}, we point out the following consequence of Theorems \ref{near} and \ref{near1}.

\begin{corollary}\label{near2}
Let $c\in\R$ be a topologically relevant critical value of a functional $J\in C^1(X,\R)$. Then,
for every $\eta >0$ there exists $\delta>0$ such that every functional $G\in C^1(X,\R)$ 
satisfying $(PS)$ in $]c-\eta, c+ \eta[$ with
\[
\sup\{|J(u)-G(u)|: u\in X\} <\delta
\]
admits a critical value in $]c-\eta, c+ \eta[$.
\end{corollary}


\section{The symmetric case}\label{s21}

In this section we deal with some existence and multiplicity results about the unperturbed problem
\[
\leqno{(P_{0})}\qquad\qquad
\left\{
\begin{array}{lll}
\displaystyle{ -\Delta u -\lambda u=  f(x, u) } & & \mbox{ in }  \Omega,\\
 & \\
 \displaystyle{u=0} & & \mbox{ on }  \partial\Omega.\\
 \end{array}
\right.
\]

In previous references on this topic, problem $(P_{0})$ has been mainly 
studied when $f=f(u)$ or $f=f(x,u)\in C(\overline\Omega\times\R, \R)$; 
instead here we deal with a Carath\'eodory nonlinearity $f$ only measurable on $x$, 
thus we need a global control on its growth. 

Remark that from $(H_0)$, $(H_1)$ and direct computations it follows that 
for all $\sigma>0$ there exists $A_\sigma>0$ such that
\begin{equation}\label{infy}
|f(x,t)|\leq\sigma|t| + A_\sigma \quad \hbox{ for all } t\in\R \hbox{ and a.e. }  x\in \Omega.
\end{equation}
Thus, the functional $I$ defined in \eqref{1} is $C^1$ (cf. e.g. \cite[Appendix C]{s}). Hence, 
the weak solutions of problem $(P_0)$ are the critical points of $I$. Let us point out that 
this is true even if, instead of \eqref{sup}, we only require that
\[
\sup_{|t|\leq {r}}|f(\cdot,t)| \in L^{\frac{2N}{N+2}}(\Omega) \quad \hbox{ for all }\quad r>0,
\]
since $\frac{2N}{N+2}$ is the conjugate exponent of $2^\ast$.
\vspace{0,5cm}

In the unperturbed non resonant case, the following result holds.

\begin{theorem}\label{multip}
Assume that $(H_0)$, $(H_1)$ and $(H_3)$ hold. Then, problem $(P_{0})$ has at least 
one solution. Moreover, under the further assumptions $(H_2)$, $(H_4)$ and $(H_5)$, problem $(P_{0})$ has at least 
$\dim (M_h\oplus\ldots \oplus M_k)$ distinct pairs of nontrivial solutions.
\end{theorem}

\dimo
Let us point out that by
$(H_3)$ functional $I$ satisfies $(PS)$ 
in $\R$ (cf. e.g. \cite{az}). Hence, a standard application of the saddle point theorem
(see \cite[Theorem 4.6]{r}) allows us to prove the existence of one solution for $(P_0)$. 
 
Here, we focus on the multiplicity statement concerning the critical points of $I$, which is even by assumption $(H_5)$. 

At first we restrict to the case $f'(0)<0$. Hence, without loss of generality, in $(H_4)$ 
we can assume  that $h$ and $k$ are such that
\begin{equation}\label{conh}
\lambda_{h-1}< f'(0) + \lambda < \lambda_h< \lambda_k <\lambda<\lambda_{k+1},
\end{equation}
possibly with $\lambda_0=-\infty$.
Notice that by \eqref{infy}, fixing any $\sigma>0$, $B_\sigma>0$ exists such that 
\[
I(u)\leq \frac 12\|u\|^2 -\frac \lambda 2|u|_2^2 + \frac \sigma 2|u|^2_2 + B_\sigma|u|_2 \qquad \hbox{for all }  u\in H^1_0(\Omega).
\]
Then, 
it results
\[
I(u)\leq \frac 12\sum_{j=1}^{k}(\lambda_j - \lambda + \sigma)|u_j|_2^2 + B_\sigma|u|_2 
\quad \hbox{ for all }  u= \sum_{j=1}^{k}u_j\in H^-(k).
\]
Hence, taking $\sigma$ small enough, by \eqref{conh} the functional $I$ tends to $-\infty$ 
as $\|u\|$ diverges in $H^-(k)$, so there exists $c_\infty\in\R$ such that $I(u)\leq c_\infty$ for all $u\in H^-(k)$.

Next we need a control on $F$ near to $t=0$.
>From $(H_1)$ and $(H_2)$, for any $\sigma>0$ there exist $R_\sigma, \delta_\sigma>0$ 
(without loss of generality $R_\sigma\geq 1$) such that
\begin{eqnarray*}
&& |F(x,t)|\leq \frac\sigma 2|t|^2 \quad \quad \hbox{if $|t|> R_\sigma$, for a.e. $x\in \Omega$, }\\
&&|F(x,t)-\frac {f'(0)}2t^2|\leq \frac\sigma 2|t|^2 \quad \quad \hbox{if $|t|< \delta_\sigma$, 
for a.e. $x\in\Omega$.}\end{eqnarray*}
Moreover, taking any $s\in [0, 2^\ast-2[$,  by \eqref{sup} there exists $a_{R_\sigma}>0$ such that,
 if $\delta_\sigma\leq |t|\leq R_\sigma$ and for a.e. $x\in\Omega$, we have
\[
|F(x,t)|\leq \sup_{|t|\leq {R_\sigma}}|f(x,t)|R_\sigma\leq a_{R_\sigma}R_\sigma
\leq a_{R_\sigma}R_\sigma\left(\frac{|t|}{\delta_\sigma}\right)^{s+2}.
\]
Summing up, for any $\sigma>0$ there exists $a_\sigma>0$ large enough such 
that for all $t\in\R$ and for a.e.  $x\in\Omega$ we have
\begin{equation}\label{conhsip}
-\frac{(\sigma - f'(0))}2|t|^2 - {a_\sigma}|t|^{s+2}\leq F(x,t)
\leq\frac{(\sigma + f'(0))}2|t|^2 + {a_\sigma}|t|^{s+2}, 
\end{equation}
which in particular implies
\[
\int_\Omega F(x,u) \; {\rm d}x\leq \frac{(\sigma + f'(0))}2|u|^2_2 +
{a_\sigma}|u|^{s+2}_{s+2} \quad \hbox{for all }  u\in H^1_0(\Omega).
\]
By the Sobolev inequalities it results
\begin{equation}\label{nicola}
I(u)\geq \frac 12\left(\|u\|^2 - (\lambda + f'(0) + \sigma)|u|_2^2\right) - 
a_\sigma' \|u\|^{s+2}\quad \hbox{for all }  u\in H^1_0(\Omega),
\end{equation}
for a suitable $a_\sigma'>0$.
On the other hand, if ${\displaystyle u= \sum_{j=h}^{+\infty}u_j\in H^+(h)}$, \eqref{nicola} implies
\[
I(u)\geq\frac 12 \sum_{j=h}^{+\infty}\left(\lambda_j - (\lambda + f'(0) + \sigma)\right)|u_j|_2^2 - a_\sigma' \|u\|^{s+2},
\]
hence, by \eqref{conh} and for $\sigma$ small enough, there exists $a''_\sigma>0$ such that 
\[
I(u)\geq  a''_\sigma\|u\|^2- a_\sigma'\|u\|^{s+2} \quad \hbox{for all }   u\in H^+(h).
\]
So, setting
$S_\rho=\{u\in H^1_0(\Omega): \|u\|= \rho\}$, if $\rho$ is small enough
there exists $c_0>0$ such that 
$I(u)\geq c_0$ for all $u\in S_\rho\cap H^+(h)$.

Moreover, considering the pseudo--index theory 
$(S_\rho\cap H^+(h), {\mathcal H}^\ast, \gamma^\ast)$ related to the genus, 
$S_\rho\cap H^+(h)$ and $I$, we have
\[
\gamma^\ast(H^-(k)) = \min _{h\in {\mathcal H}^\ast}\gamma(H^-(k)\cap h^{-1}(S_\rho\cap H^+(h)))
\geq \dim H^-(k) - \mbox{ codim }H^+(h)
\]
(cf. \cite[Theorem A.2]{bbf}).
Hence, Theorem \ref{group} applies with $\bar A:= H^-(k)$ and $S:=S_\rho\cap H^+(h)$, 
so $I$ has at least $\dim (M_h\oplus\ldots \oplus M_k)$ distinct pairs of critical 
points corresponding to at most $\dim (M_h\oplus\ldots \oplus M_k)$ distinct 
critical values $c_i$, where $c_i$ is as in \eqref{value}.

If $f'(0)>0$, using \eqref{conhsip}, the proof follows by applying 
Theorem \ref{group} to the functional $-I$ with $\bar A=:H^+(h)$ and $S:=S_\rho\cap H^-(k)$.
Let us point out that in this case $H^+(h)$ is infinite dimensional, then by \cite[Theorem 3.4]{bcf} it is
\begin{eqnarray*}
\gamma^\ast(H^+(h)) &\geq& \dim(H^-(k)\cap  H^+(h)) - \mbox{ codim } (H^-(k) + H^+(h)) \\
& = & \dim (M_h\oplus\ldots \oplus M_k), 
\end{eqnarray*}
thus the proof is complete. \hfill$\square$

\begin{remark}\label{inzero}
{\em If in assumption $(H_2)$ it is $f'(0)=0$, replacing assumption $(H_4)$ by
\[
F(x,t) < 0 \quad \hbox{ for a.e. } \quad x\in\Omega  \quad \hbox{ and } \quad t\not= 0,
\]
and reasoning as in \cite[Theorem 6.1]{bbf}, it results $I(u)\geq c_0$ for 
all $u\in S_\rho\cap H^+(k)$. Hence we obtain $\dim H^-(k) -\mbox{ codim } H^+(k) = \dim M_k$ distinct pairs of critical points.}
\end{remark}

\begin{remark}\label{dora2}
{\em The multiplicity result stated in Theorem \ref{multip} still holds if the limit 
in assumption $(H_2)$ is infinite. More precisely, if the assumptions $(H_0)$,  
$(H_1)$, $(H_3)$, $(H_5)$ hold and moreover we assume
\begin{itemize}
\item[$(H_7)$]
there exists
\[
\lim_{t\rightarrow 0}\frac{f(x,t)}{t} \ =-\infty \quad \hbox{uniformly for a.e. 
$x \in \Omega$ and $\lambda_k<\lambda$ for some $k\in\N$,}\]
\end{itemize}
then   $(P_0)$ has at at least $\dim (M_1\oplus\ldots \oplus M_k)$ distinct 
pairs of non--trivial solutions, choosing in Theorem \ref{group} 
$\bar A = H^-(k)$ and $S=S_\rho\cap H^+(1)$, i.e. $S=S_\rho$.

On the other hand, if $(H_7)$ above is replaced by
\begin{itemize}
\item[$(H_7')$]
there exists
\[
\lim_{t\rightarrow 0}\frac{f(x,t)}{t} \ =+\infty \quad \hbox{uniformly for a.e. 
$x \in \Omega$,}\]
\end{itemize}
then $(P_0)$ has infinitely many pairs of non--trivial solutions. 
Indeed, fixing $h$ such that $\lambda_h>\lambda$, for any $k>h$ we 
can apply Theorem \ref{group} to the functional $-I$ with 
$\bar A=H^+(h)$ and $S=S_\rho\cap H^-(k)$, thus obtaining 
$\dim (M_h\oplus\ldots \oplus M_k)$ pairs of solutions. The conclusion follows by the arbitrariness of $k$.}
\end{remark}

\begin{remark}
{\em Theorem \ref{multip} holds also in the resonant case, with assumption $(H_1)$ 
replaced by the stronger one $(H_1')$ and   $(H_3)$, $(H_4)$ replaced respectively 
by $(H_3')$, $(H_4')$, but adding $(H_6)$ (see also \cite[Theorem 4.12]{r}). 
Moreover, the arguments in Remarks \ref{inzero} and \ref{dora2} still work. }
\end{remark}

\section{Proof of the main results}\label{s23}

\noindent
{\sl Proof of Theorem \ref{maintheorem}.\ }
Here, we prove only the multiplicity result, since simpler arguments give 
the existence of one solution. Indeed, by the first statement in Theorem  
\ref{multip} we find a critical point of the unperturbed functional $I$ 
in \eqref{1} and by Theorem \ref{closed} the corresponding critical value 
is  an essential one of $I$, hence by Theorems \ref{near} and \ref{near1} the thesis follows.

Assume that $f'(0)<0$ and \eqref{conh} holds (similar arguments work if $f'(0)>0$).
Fixing any $j\in\N$, as in \cite{ll} we consider a continuous cut function 
\[
\beta_j\left(t\right)=
\left\{
\begin{array}{lll}
\displaystyle{0 } & {\rm if} & |t|\geq j+1\\
 & \\
 \displaystyle{1} & {\rm if} & |t|\leq j\\
 \end{array}
\right.
\]
such that $0<\beta_j(t)<1$ if $j< |t| < j+1$. Then, let us set
\[
g_j(x,t) =\beta_j(t)g(x,t) \qquad \hbox{and} \qquad G_j(x,t) = \int_0^t g_j(x,s)\; {\rm d}s.
\]
Let us remark that, if $g(x,\cdot)$ is odd, then choosing $\beta_j$ even, 
it results that $g_j(x,\cdot)$ and $G_j(x,\cdot)$ are odd and even respectively, for a.e. $x\in\Omega$.
Furthermore, there exists $\varepsilon_1(j)>0$ such that 
\begin{equation}
\label{2}  
\varepsilon_1(j)|g_j(x,t)|<1, \qquad \varepsilon_1(j)|G_j(x,t)|<1 \quad \hbox{ for all } x\in\Omega, t\in\R,
\end{equation}
thus for any $|\eps|\leq \varepsilon_1(j)$ we consider the functionals
\[
I_{j,\varepsilon}(u)= I(u) - \varepsilon\int_\Omega  G_j(x,u)\; {\rm d}x \quad \hbox{ on } H^1_0(\Omega).
\]

Let $\bar m$ be the number of the distinct critical levels $c_i$ of $I$ found 
in Theorem \ref{multip}. Clearly $1\leq\bar m\leq \dim (M_h\oplus\ldots \oplus M_k)$ 
and $0<c_0<c_{i_1}<\ldots<c_{i_{\bar m}}\leq c_\infty$, where $c_0$ and $c_\infty$ are 
as in the proof of Theorem \ref{multip}.
These critical levels are also odd--essential ones for $I$.
Indeed, in order to prove this it suffices to apply 
Corollary \ref{near5}. Namely, we take $X=H^1_0(\Omega)$, 
$\Gamma=\Sigma_i$ defined in Theorem \ref{group}, $d=0$. 
Then, for any odd homeomorphism $\varphi:H^1_0(\Omega)\times [0,1]\longrightarrow H^1_0(\Omega)$ 
such that $\varphi(u,t)=u$ if $I(u)\leq 0$, we have that the set $\overline{\varphi(C \times\{1\})}$
is closed and symmetric,  for each $C\in\Sigma_i$. Moreover, from the supervariancy property of $\gamma^*$, it is:
\[
\gamma^\ast\left(\overline{\varphi(C \times\{1\})}\right)=
\gamma^\ast\left(\varphi\left(\overline{C \times\{1\}}\right)\right)=
\gamma^\ast\left(\overline{C \times\{1\}}\right)=\gamma^\ast (C)\geq i,
\]
hence $\overline{\varphi(C \times\{1\})}$ belongs to $\Sigma_i$ and the conclusion follows.

So, by Theorem \ref{near} there exists $\varepsilon_2(j)\in ]0, \varepsilon_1(j)[$ such 
that, if $|\varepsilon|\leq \varepsilon_2(j)$, then
$I_{j,\varepsilon}$ has at least $\bar m$ odd--essential values $d_i^{j,\varepsilon}$, with $i\in \{1,\ldots, \bar m\}$, such that
\begin{equation}
\label{3}
\frac{c_0}2< d_1^{j,\varepsilon}< \ldots < d_{\bar m}^{j,\varepsilon}< c_\infty +1.
\end{equation}
As $I_{j,\varepsilon}$ satisfies the Palais--Smale condition in $\R$ (cf. e.g. \cite{az}),
by Theorem \ref{near1},
for each  $i\in \{1,\ldots, \bar m\}$ $I_{j,\varepsilon}$ has a critical point $u_i^{j,\varepsilon}$ such that 
\[
(P_{j,\varepsilon})\qquad\qquad
\left\{
\begin{array}{lll}
\displaystyle{ -\Delta u_i^{j,\varepsilon} -\lambda u_i^{j,\varepsilon}=  f(x, u_i^{j,\varepsilon}) 
+ \varepsilon g_j(x, u_i^{j,\varepsilon})  } & & \mbox{ in }  \Omega,\\
 \displaystyle{u_i^{j,\varepsilon}=0} & & \mbox{ on }  \partial\Omega\\
 \end{array}
\right.
\]
and
\[
d^{j,\varepsilon}_i =\frac 12\int_\Omega |\nabla u_i^{j,\varepsilon}|^2\; {\rm d}x - 
\frac\lambda 2 \int_\Omega (u_i^{j,\varepsilon})^2 \; {\rm d}x - \int_\Omega \left(F(x,u_i^{j,\varepsilon}) + \varepsilon G_j(x,u_i^{j,\varepsilon})\right) \; {\rm d}x.
\]

We claim that 
\begin{equation}
\label{6}
\| u_i^{j,\varepsilon}\|\leq C_1 \quad \hbox{ for all } j\in \N,|\varepsilon|\leq \varepsilon_2(j), i\in \{1,\ldots, \bar m\}.
\end{equation}
Firstly, for all $u\in H^1_0(\Omega)$, as $u=u^+ + u^-$ with $u^+\in H^+(k+1)$ and $u^-\in H^-(k)$, 
standard computations show that there exists $\delta>0$ such that
\begin{eqnarray}
\label{4.6} 
&&   \|u^+\|^2 - \lambda|u^+|_2^2= \sum_{i=k+1}^{+\infty}(\lambda_i - \lambda)|u_i|_2^2\geq \delta\|u^+\|^2,            \\
\label{4.7}
&&   \lambda|u^-|^2 -\|u^-\|_2^2= \sum_{i=1}^{k}(\lambda - \lambda_i)|u_i|_2^2\geq \delta\|u^-\|^2.\end{eqnarray}
Clearly, from $(P_{j,\varepsilon})$ and \eqref{2} we have
\begin{eqnarray*}
&& 
\|(u_i^{j,\varepsilon})^+\|^2-\lambda |(u_i^{j,\varepsilon})^+|_2^2\leq \int_\Omega
|f(x,u_i^{j,\varepsilon})||(u_i^{j,\varepsilon})^+|\; {\rm d}x +|(u_i^{j,\varepsilon})^+|_1, \\
&&
\lambda |(u_i^{j,\varepsilon})^-|_2^2 - \|(u_i^{j,\varepsilon})^-\|^2\leq \int_\Omega
|f(x,u_i^{j,\varepsilon})||(u_i^{j,\varepsilon})^-|\; {\rm d}x +|(u_i^{j,\varepsilon})^-|_1.\end{eqnarray*}
Hence, by \eqref{infy} and \eqref{4.6}, respectively \eqref{4.7},  for suitable $\tilde\eps, C_2>0$, by standard computations we obtain
\[
(\delta- \tilde\eps)\|(u_i^{j,\varepsilon})^+\|^2\leq \tilde\eps \|u_i^{j,\varepsilon}\|^2 + C_2\|u_i^{j,\varepsilon}\|,
\]
respectively
\[
(\delta- \tilde\eps)\|(u_i^{j,\varepsilon})^-\|^2\leq \tilde\eps \|u_i^{j,\varepsilon}\|^2 + C_2
\|u_i^{j,\varepsilon}\|.
\]
Putting together the two previous inequalities and choosing $\tilde\eps$ small enough, we get that \eqref{6} holds.

Using regularity results and a standard bootstrap method, we show that from \eqref{6} we have
\begin{equation}
\label{7}
| u_i^{j,\varepsilon}|_\infty\leq C_3 \quad\hbox{ for all } j\in \N, |\varepsilon|\leq\varepsilon_2(j),i\in \{1,\ldots, \bar m\}.
\end{equation}
Indeed, for $j$ fixed, $u_i^{j,\varepsilon}\in L^{2^\ast}(\Omega)$ and so, as $f$ is sublinear, also
\[
\phi^{j,\varepsilon}_i(x,u_i^{j,\varepsilon}):= \lambda u_i^{j,\varepsilon} + f(x,u_i^{j,\varepsilon}) 
+\varepsilon g_{j}(x, u_i^{j,\varepsilon})
\]
belongs to $L^{2^\ast}(\Omega)$.
Then, by \cite[Theorem B.2]{s}, it follows that $u_i^{j,\varepsilon}\in 
H_0^{2,2^\ast}(\Omega)$ and
\[
\| u_i^{j,\varepsilon}\|_{H_0^{2,2^\ast}}\leq C_4\left(| u_i^{j,\varepsilon}|_{{2^\ast}} + | \phi^{j,\varepsilon }_i(x,u_i^{j,\varepsilon})|_{{2^\ast}}\right).
\]
Then from \eqref{6} we get 
\[
\| u_i^{j,\varepsilon}\|_{H_0^{2,2^\ast}}\leq C_5.
\]
Now, from \cite[Theorem A.5]{s}, if $N\leq 5$, \eqref{7} is true. Otherwise, if $N>6$ ($N=6$ is a simpler case), as
$H_0^{2,2^\ast}(\Omega)$ is continuously embedded in $L^{q^\ast}(\Omega)$, with $q^\ast=\frac {2N}{N-6}$, we have that
$u_i^{j,\varepsilon}$ and $\phi^{j,\varepsilon}_i(x,u_i^{j,\varepsilon})$ are in $H_0^{2,q^\ast}(\Omega)$ and, applying again 
\cite[Theorem B.2]{s}, we get \eqref{7} for $N\leq 9$. Going on in this way, \eqref{7} holds for any $N\in\N$.
Finally, for $j> C_3$, problem $(P_{\varepsilon})$ has at least $\bar m$ pairs of solutions. \hfill$\square$

\begin{remark}\label{dora3}
{\em By Remarks \ref{inzero}, \ref{dora2} and corresponding suitable changes, 
the result in Theorem \ref{maintheorem} still holds if in $(H_2)$ it is 
$f'(0)=0$ or $f'(0)\in\{\pm \infty\}$, with $\bar m\geq 1$ number of the distinct critical levels of the unpertubed functional.}
\end{remark}

\noindent
{\sl  Proof of Theorem \ref{strong}.\ }
Let us consider the cut functions $\beta_j$ and the notations as in the proof of Theorem \ref{maintheorem}.
Under our assumptions the functionals $I$ and $I_{j,\varepsilon}$ satisfy $(PS)$ 
(see \cite[Theorem 4.12]{r}).  Then $I_{j,\varepsilon}$ has at least $\bar m$ 
odd--essential values $d^{j,\varepsilon}_i$ verifying \eqref{3}.
Now each $u\in H^1_0(\Omega)$ can be written as $u=u^+ + u^- + u^0$, with $u^+\in H^+(k+1), u^-\in H^-(k-1)$ and $u_0\in M_k$.
Again standard computations show that there exists $\delta>0$ such that
\begin{eqnarray}
\label{4.6r} 
&&   \|u^+\|^2 - \lambda_k|u^+|_2^2\geq \delta\|u^+\|^2,            \\
\label{4.7r}
&&   \lambda_k|u^-|^2 -\|u^-\|_2^2\geq \delta\|u^-\|^2.\end{eqnarray}
>From $(H_1')$, \eqref{2}, $(P_{j,\varepsilon})$ and \eqref{4.6r}--\eqref{4.7r}, it follows that 
\[
\|(u_i^{j,\eps})^\pm\|
\leq C_1 \quad\hbox{ for all } j\in \N, |\varepsilon|\leq\varepsilon_2(j),i\in \{1,\ldots, \bar m\}.
\]
We claim that also $\|(u_i^{j,\eps})^0\|$ is bounded. Indeed, again $(H_1')$ and standard arguments imply
\[
\left|\int_\Omega F(x,u)\; {\rm d}x\right|\leq C_2\|u\| \quad \hbox{ for all } \quad u \in H^1_0(\Omega);
\]
thus, as
\[
I_{j,\eps}(u_i^{j,\varepsilon})\leq c_\infty + 1  \quad\hbox{ for all } j\in \N, 
|\varepsilon|\leq \varepsilon_2(j), i\in \{1,\ldots, \bar m\},
\]
the thesis follows by \eqref{2} and reasoning as in \cite[Theorem 4.12 and Lemma 4.21]{r}.
\hfill$\square$

\vspace{0,5cm}
\noindent
{\sl Proof of Theorem \ref{dis}.\ }
In this case, once found the critical values of $I$ by Theorem \ref{multip}, 
as they are assumed to be topologically relevant, we can apply Corollary \ref{near2}:
so there exists $\varepsilon_2(j)\in ]0, \varepsilon_1(j)]$ such that, if $|\varepsilon|\leq \varepsilon_2(j)$, then
$I_{j,\varepsilon}$ has at least $\bar m$ critical values $d_k^{j,\varepsilon}$, with $k\in \{1,\ldots, \bar m\}$, such that
\eqref{3} holds. Then we proceed as in the proof of Theorem \ref{maintheorem}.\hfill$\square$

\vspace{0,5cm}
\noindent
{\sl Proof of Theorem \ref{disst}.\ }
It is enough to combine the arguments in the proofs of Theorems \ref{strong} and \ref{dis}.\hfill$\square$

\end{document}